\documentclass[11pt,letterpaper,reqno]{amsart}

\usepackage{amssymb,amsfonts}
\usepackage[all,arc]{xy}
\usepackage{enumerate}
\usepackage{mathrsfs}
\usepackage{graphicx}
\usepackage{cite}

\newtheorem{theorem}{Theorem}[section]
\newtheorem{corollary}[theorem]{Corollary}
\newtheorem{proposition}[theorem]{Proposition}
\newtheorem{lemma}[theorem]{Lemma}

\theoremstyle{definition}
\newtheorem{definition}[theorem]{Definition}

\newtheorem{question}[theorem]{Question}
\newtheorem{remark}[theorem]{Remark}

\let\oldbigwedge\bigwedge
\def\BIGwedge{{\textstyle\oldbigwedge}}
\def\medwedge{{\scriptstyle\oldbigwedge}}
\def\bigwedge{\mathchoice{\BIGwedge}{\BIGwedge}{\medwedge}{}}

\DeclareMathOperator{\diam}{diam}

\DeclareMathOperator{\Ann}{Ann}
\DeclareMathOperator{\id}{id}
\DeclareMathOperator{\Id}{Id}
\DeclareMathOperator{\Sub}{Sub}


\numberwithin{equation}{section}

\begin{document}
	
	\title{A generalization of zero-divisor graphs}
	
	\author{Peyman Nasehpour}
	
	\address{Peyman Nasehpour, Department of Engineering Science, Golpayegan University of Technology, Golpayegan, Iran}
	
	\email{nasehpour@gut.ac.ir, nasehpour@gmail.com}

	

	\begin{abstract}
		In this paper, we introduce a family of graphs which is a generalization of zero-divisor graphs and compute an upper-bound for the diameter of such graphs. We also investigate their cycles and cores.
	\end{abstract}

	\maketitle
	
	\section{Introduction}\label{sec:intro}
	
	For coloring a commutative ring, Beck introduced a version of the zero-divisor graph of a ring in his 1988 paper \cite{Beck1988}. Later in 1999, Anderson and Livingston introduced a similar notion which is the by-now standard definition of zero-divisor graphs \cite{AndersonLivingston1999}. This notion has been generalized and investigated for commutative semigroups with zero by DeMeyer et al.  \cite{DeMeyerMcKenzieSchneider2002,DeMeyerDeMeyer2005}. Since then, many authors have investigated the zero-divisor graphs from different perspectives and for a survey on this, one may refer to the papers \cite{AndersonAxtellStickles2011,AndersonBadawi2017}. Similarly, for non-commutative rings, Redmond has introduced a similar notion called zero-divisor (directed) graphs \cite{Redmond2002}. 
	
	One of the interesting topics in algebraic combinatorics is to compute invariants of zero-divisor graphs such as their diameters, girths, clique numbers, chromatic numbers, and even ``Zagreb indices'' \cite{AykacAkgunesCevik2019} and for a survey on the computation of these invariants, one can check the paper \cite{CoykendallWagstaffSheppardsonSpiroff2012}. For the comparison of these numbers for zero-divisor graphs of a semigroup under Armendariz extension one may see the 2013 paper by Epstein et al. \cite{EpsteinNasehpour2013} and under polynomial and power series extensions the 2006 paper by Lucas \cite{Lucas2006}. Section 5 of the 2010 paper \cite{Nasehpour2010} is devoted to the comparison of the diameter of zero-divisor graphs under content extensions. One interesting topic for a future project can be to compute the tenacity \cite{CozzensMoazzamiStueckle1995} of zero-divisor graphs.
	
	Our main motivation for this paper was to attribute a graph $RG(M)$ to a module $M$ inspired by zero-divisor graphs of ideals of a ring in the following sense:
	
	Let $R$ be a commutative ring with a nonzero identity and $M$ be a unital $R$-module. We associate a graph $RG(M)$ to $M$, which we call residuated graph of $M$, whose vertices and edges are determined as follows:
	
	\begin{enumerate}
		\item Let $N$ be a submodule of $M$. Then $N$ is a vertex of $RG(M)$ if the residuated ideal $[N:M]$ of $R$ is nonzero and there is a submodule $K\neq N$ of $M$ with $[K:M] \neq (0)$ such that \[[N:M]\cdot [K:M]=(0),\] where by $[N:M]$, we mean the set of all elements $r\in R$ such that $rM \subseteq N$;
		
		\item Two distinct vertices $P$ and $Q$ of the graph $RG(M)$ are connected if \[[P:M]\cdot [Q:M]=(0).\] 
	\end{enumerate}
	
	Surprisingly, similar to the zero-divisor graphs of commutative semigroups \cite[Theorem 1.3]{DeMeyerMcKenzieSchneider2002}, the graph $RG(M)$, for any $R$-module $M$, is connected and the best upper-bound for $\diam RG(M)$ is 3 if the graph $RG(M)$ is non-empty (see Corollary \ref{diamsemimodule}). Here we need to recall that the distance between two vertices in a simple graph is the number of edges in a shortest path connecting them. The greatest distance between any two vertices in a graph $G$ is the diameter of $G$, denoted by $\diam(G)$ \cite[p. 8]{Diestel2017}.
	
	Based on our investigations for residuated graphs, in Definition \ref{zdgDef}, we attribute a graph to an arbitrary set which is also a generalization of the notion of zero-divisor graphs of arbitrary commutative semigroups with zero in the following sense:
	
	Let $X$ be a non-empty set, $(S, \cdot, 0)$ a commutative multiplicative semigroup with zero, and $f$ a function from $X$ to $S$. We attribute a simple graph to $X$, denoted by $\Gamma_{(S,f)}(X)$, whose vertices and edges are determined as follows:
	
	\begin{enumerate}
		
		\item An element $x\in X$ is a vertex of the graph $\Gamma_{(S,f)}(X)$ if $f(x) \neq 0$ and there is a $y\neq x$ in $X$ such that $f(y) \neq 0$ and $f(x)\cdot f(y) =0$.
		
		\item Let $x$ and $y$ be elements of $X$. The doubleton $\{x,y\}$ is an edge of the graph $\Gamma_{(S,f)}(X)$ if $x \neq y$, $f(x) \neq 0$, and $f(y)\neq 0$ while $f(x)\cdot f(y) = 0$.
	\end{enumerate}
	
	Then, in Section \ref{sec:zdg}, we prove that under some conditions, the graph $\Gamma_{(S,f)}(X)$ is connected with $\diam \Gamma_{(S,f)}(X) \leq 3$ if $\Gamma_{(S,f)}(X)$ is non-empty (see Definition \ref{zdgDef}, Theorem \ref{diamsemigroup}, and Theorem \ref{diamsemigroup2}).

	Note that in the Definition \ref{zdgDef}, if we set $X=S$ and suppose that $\id_S$ is the identity map on a commutative semigroup with zero $S$, then $\Gamma_{(S,\id_S)}(S)$ is nothing but the zero-divisor graph $\Gamma(S)$ defined in \cite{DeMeyerMcKenzieSchneider2002}.
	
	In Section \ref{sec:diam}, we prove that if $S$ is a commutative semiring with a nonzero identity and the $S$-semimodule $M$ has the annihilator condition or $M$ is a content $S$-semimmodule and the content function from $M$ to finitely generated ideals of $S$ is onto, then the graphs $\Gamma_{(\Id(S),\Ann)}(M) $ and $\Gamma_{(\Id(S), c)}(M)$ are connected with diameters at most 3 if they are non-empty (see Corollary \ref{zdgac} and Corollary \ref{zdgcontent}).
		
    We also show that if $S$ is a commutative semiring with a nonzero identity, $M$ is a unital $S$-semimodule, $q$ is a function from $\Sub(M)$ to $\Id(S)$ with $q(N) = [N:M]$, and the graph $\Gamma_{(\Id(S),q)}(\Sub(M))$ is non-empty, then it is a connected graph whose diameter is at most 3 (see Corollary \ref{diamsemimodule}).
    
    In Section \ref{sec:core}, we discuss the cycles and cores of the graphs defined in Definition \ref{zdgDef}. For example in Theorem \ref{core}, we prove that if $X$ is a non-empty set, $S$ a commutative semigroup with zero, $f$ a function from $X$ to $S$, the graph $\Gamma_{(S,f)}(X)$ has at least three vertices, and the function $f$ has this property that for all $x,y\in X$ if $f(x)f(y)\neq 0$ then there exists a $z\in X$ such that $f(z)=f(x)f(y)$, then if $\Gamma_{(S,f)}(X)$ contains a cycle, then the core $K$ of $\Gamma_{(S,f)}(X)$ is a union of triangles and rectangles.
	
	We recall that a trail in a graph $G$ is a walk in which all edges are distinct. A path in the graph $G$ is a trail in which all vertices (except possibly the first and last) are distinct. If $ P = x_0 \cdots x_{k-1}$ is a path in $G$ and $k \geq 3$, then the path $C = x_0 \cdots x_{k-1} x_0$ is a cycle in $G$ \cite{Diestel2017}. We also note that the core of a graph $\Gamma$ is the largest subgraph of $\Gamma$ in which every edge is the edge of a cycle in $\Gamma$ \cite{DeMeyerDeMeyer2005}.
	
	\section{A generalization of zero-divisor graphs for semigroups}\label{sec:zdg}
	
	One of the interesting areas of research in algebraic combinatorics is to associate a graph $G(A)$ to an algebraic structure $A$ and investigate the interplay between the algebraic properties of the algebra $A$ and the graph-theoretic properties of the graph $G(A)$.  One method is to consider the intersection graphs of the substructures of an algebraic structure. For example, in the 2012 paper \cite{AkbariTavallaeeGhezelahmad2012}, Akbari et al. investigate the intersection graphs of the submodules of modules over arbitrary commutative rings. Since 1960s, many authors have worked on intersection graphs \cite{Bosak1964, ChakrabartyGhoshMukherjeeSen2009, CsakanyPollak1969, Osba2016, Shen2010, Zelinka1975, Zelinka1973}. Note that all graphs are intersection graphs \cite{ErdosGoodmanPosa1966}. In this direction, Malakooti Rad and Nasehpour generalize the notion of intersection graphs and attribute a graph to the bounded semilattices and investigate their properties and compute the invariants of such graphs \cite{MalakootiRadNasehpour2017}.
	
	In this section, we attribute a graph to an arbitrary set which is on one hand a generalization of the notion of zero-divisor graphs of commutative semigroups and on the other hand is a generalization of the graphs attributed to submodules of a module given in Corollary \ref{diamsemimodule}.
	
	\begin{definition}
		
		\label{zdgDef}
		
		Let $X$ be a non-empty set, $(S, \cdot, 0)$ a commutative multiplicative semigroup with zero, and $f$ a function from $X$ to $S$. We attribute a graph to $X$, denoted by $\Gamma_{(S,f)}(X)$, whose vertices and edges are determined as follows:
		
		\begin{enumerate}
			\item An element $x\in X$ is a vertex of the graph $\Gamma_{(S,f)}(X)$ if $f(x) \neq 0$ and there is a $y\neq x$ in $X$ such that $f(y) \neq 0$ and $f(x)\cdot f(y) =0$.
			
			\item Let $x$ and $y$ be elements of $X$. The doubleton $\{x,y\}$ is an edge of the graph $\Gamma_{(S,f)}(X)$ if $x \neq y$, $f(x) \neq 0$, and $f(y)\neq 0$ while $f(x)\cdot f(y) = 0$.
		\end{enumerate}
		
	\end{definition}
	
	\begin{remark}
			
	Let $X$ be a non-empty set, $S$ a commutative semigroup with zero, and $f$ a function from $X$ to $S$. The graph $\Gamma_{(S,f)}(X)$ is a generalization of the usual zero-divisor graph $\Gamma(S)$ defined in \cite{DeMeyerMcKenzieSchneider2002}. In fact, if suppose that $S$ is a commutative semigroup with zero and $X=S$, then $\Gamma_{(S,\id_S)}(S)$ is the zero-divisor graph $\Gamma(S)$, where $\id_S$ is the identity map on $S$.
	
	A graph $C$ is called to be a zero-divisor if these exist non-isomorphic graphs $A$ and $B$ for which $A \times C \cong B \times C$ \cite[p. 310]{ImrichKlavzar2000}. For examples of these graphs see \cite{Lovasz1971}. And one should not confuse this concept in graph theory with the concept of zero-divisor graphs in \cite{DeMeyerMcKenzieSchneider2002}.  
		
	\end{remark}

\begin{question}
	Let $G$ be an arbitrary graph. Is it possible to find a set $X$, a commutative semigroup with zero $S$, and a function $f$ from $X$ to $S$ such that $G$ is isomorphic to the graph $\Gamma_{(S,f)}(X)$? 
\end{question}

\section{Diameter of Zero-Divisor Graphs and Their Generalizations}\label{sec:diam}
	
	\begin{theorem}
		
		\label{diamsemigroup}
		
		Let $X$ be a non-empty set, $S$ a commutative semigroup with zero, and $f$ a function from $X$ to $S$ with this property that for all $x,y\in X$, if $f(x)f(y)\neq 0$ then there exists a $z\in X$ such that $f(z)=f(x)f(y)$. Then the graph $\Gamma_{(S,f)}(X)$ is connected with $\diam(\Gamma_{(S,f)}(X)) \leq 3$.
		
		\begin{proof}
			Let $x,y$ be two distinct vertices of $\Gamma_{(S,f)}(X)$. Therefore, there exists $z,w \in X$ such that $f(z)\neq 0$, $f(w)\neq 0$ and $f(x)f(z)=0$ and $f(y)f(w)=0$. Note that by definition, $f(x)\neq 0$ and $f(y)\neq 0$. 
			
			Now we show that $d(x,y)\leq 3$. If $f(x)f(y)=0$, then $d(x,y)=1$. If $f(x)f(y)\neq 0$, but $f(z)f(w)=0$, then $x-z-w-y$ is a path in $\Gamma_{(S,f)}(X)$ and therefore, $d(x,y)\leq 3$.
			
			Finally, let $f(x)f(y)\neq 0$ and $f(z)f(w)\neq 0$. Since there exists a $t \in X$ such that $f(t)=f(z)f(w)$, we have $f(x)f(t) = f(t)f(y) = 0$ and $d(x,y)\leq 2$. Therefore, the graph $\Gamma_{(S,f)}(X)$ is connected with diameter at most 3 and the proof is complete.
		\end{proof}
		
	\end{theorem}

\begin{corollary}
	
	\label{diamsemigroupcor}
	
	Let $S$ be a commutative semigroup with zero. The zero-divisor graph $\Gamma(S)$ is connected with $\diam \Gamma(S) \leq 3$ \cite[Theorem 1]{DeMeyerDeMeyer2005}.
\end{corollary}

Let $X$ be a non-empty set, $S$ a commutative semigroup with zero, and $f$ a function from $X$ to $S$. We do not know if the graph $\Gamma_{(S,f)}(X)$ is connected in general. Based on this, the following question arises:

\begin{question}
	Let $X$ be a non-empty set, $S$ a commutative semigroup with zero, and $f$ a function from $X$ to $S$. If the graph $\Gamma_{(S,f)}(X)$ defined in Definition \ref{zdgDef} is connected, what is the best upper-bound for the diameter of this graph?
\end{question}

Related to the above question, we bring the following remark:

\begin{remark}
	Let us recall that if $S$ is a semigroup (not necessarily commutative) with zero, a directed graph $\Gamma(S)$, called zero-divisor graph of $S$, is attributed to $S$ whose vertices are the proper zero-divisors of $S$ and $s \rightarrow t$ is an edge of $\Gamma(S)$ between the vertices $s$ and $t$ if $s\not=t$ and $st=0$ \cite{CannonNeuerburgRedmond2005}. The following result from \cite{CannonNeuerburgRedmond2005,Redmond2002}, is an interesting generalization of Corollary \ref{diamsemigroupcor} though written in the terminology of the paper \cite{Nasehpour2021}:
	
	\begin{theorem} Let $S$ be a semigroup with zero. The directed graph $\Gamma(S)$ is connected if and only if $S$ is eversible. Moreover, if $\Gamma(S)$ is connected, then the diameter of the graph $\Gamma(S)$ is at most 3.
	\end{theorem}

Note that a semigroup with zero $S$ is eversible if every left zero-divisor on $S$ is also a right zero-divisor on $S$ and conversely, i.e., $Z_l(S) = Z_r(S)$ \cite[Definition 1.9]{Nasehpour2021}.
\end{remark}
	
	Let us recall that a commutative ring $R$ with an identity has the annihilator condition if for all $a,b \in R$, there is a $c\in R$ such that $\Ann(a,b) = \Ann(c)$ \cite{HuckabaKeller1979}. Inspired by this, we give the following definition for semimodules \cite[Chap. 14]{Golan1999}:
	
	\begin{definition}
		
		\label{ac}
		
		Let $S$ be a commutative semiring with an identity and $M$ be a unital $R$-semimodule. We say that $M$ has the annihilator condition if for all $x,y \in M$, there is a $z\in M$ such that $\Ann(x,y) = \Ann(z)$, where by $\Ann(N)$, we mean the set of all elements $s$ in $S$ such that $sN=0$.
	\end{definition}

     Note that we gather all ideals of a semiring $S$ in the set $\Id_S(S)$ and all $S$-subsemimodules of $M$ in the set $\Sub_S(M)$.
	
	\begin{corollary}
		
		\label{zdgac}
		Let the $S$-semimodule $M$ have the annihilator condition. Then the graph $\Gamma_{(\Id(S),\Ann)}(M) $ is a connected graph whose dimater is at most 3.
		
		\begin{proof}
			It is clear that $(\Id(S), \cap)$ is a commutative semigroup and its zero, i.e., its absorbing element, is the zero ideal $(0)$. Consider the function $\Ann$ from $M$ to $\Id(S)$. It is straightforward to see that $\Ann(x,y)= \Ann(x) \cap \Ann(y)$ for all $x,y \in M$. Since by assumption the $S$-semimodule $M$ has the annihilator condition, the proof is complete. 
		\end{proof}
	\end{corollary}
	
	Let $S$ be a commutative semigroup with zero. A subset $I$ of $S$ is said to be an $s$-ideal of $S$, if $0\in I$ and for all $s\in S$ and $a\in I$, we have $s \cdot a \in I$ \cite{Aubert1953}. Clearly, the intersection of two $s$-ideals of a semigroup $S$ is an $s$-ideal of $S$. If we denote the set of all $s$-ideals of $S$ by $\Id_S(S)$, then $\Id_S(S)$ along with the intersection configures a commutative semigroup with zero and its absorbing element is the $s$-ideal $\{0\}$.
	
	Let us recall that if $S$ is a semigroup, a set $M$ together with a function $S \times M \rightarrow M$, denoted $(s,m) \rightarrow sm$, satisfying $(st)x = s(tx)$ for all $s,t \in S$ and $x\in M$ is called a (left) $S$-act. Also, if $M$ is a $S$-act and the semigroup $S$ has an absorbing element $0_S$ and $M$ possesses a distinguished element $0_M$ such that $s0_M = 0_M$ for all $s\in S$ and $0_S x=0_M$ for all $x\in M$, then $M$ is called a pointed $S$-act. Finally, if $S$ is a monoid and $1_S$ is the neural element for the multiplication of $S$, then an $S$-act $M$ is called unital if $1_S m = m$ for all $m\in M$ \cite{Talwar1995}. Note that if $S$ is a semiring and $M$ is a unital $S$-semimodule, then obviously $M$ is a unital pointed $S$-act.
	
	Now, let $S$ be a commutative monoid with zero and $M$ a unital pointed $S$-act. If $\emptyset \neq N \subseteq M$, we define $\Ann(N)$ to be the set of all elements $s\in S$ such that $sN=\{0_M\}$. One can easily check that $\Ann(N)$ is an $s$-ideal of the semigroup $S$ and if $P$ and $Q$ are non-empty subsets of $M$, then $\Ann(P) \cap \Ann(Q) = \Ann(P \cup Q)$. Therefore, we have already showed that the following result is just another example for Theorem \ref{diamsemigroup}:
	
	\begin{corollary}
		Let $S$ be a commutative monoid with zero and $M$ a unital pointed $S$-act. If $\mathcal{C}$ is a non-empty class of non-empty subsets of the set $M$ and $(\mathcal{C}, \cup)$ is a semigroup and the graph $\Gamma_{(\Id_S(S), \Ann)}(\mathcal{C})$ is non-empty, then it is a connected graph with diameter at most 3.  
	\end{corollary}
	
	Let us recall that if $S$ is a commutative semiring with a nonzero identity and $M$ is a unital $S$-semimodule, then the content function from $M$ into the ideals $\Id(S)$ of $S$ is defined as follows: $$c(x) = \bigcap \{I\in \Id(S): x\in IM\}.$$ An $S$-semimodule $M$ is called a content semimodule if $x\in c(x)M$ for all $x\in M$. It is straightforward to see that if $M$ is a content $S$-semimodule, then $c(x)$ is a finitely generated ideal of $S$ for each $x\in M$ \cite[Proposition 23]{Nasehpour2016}. Now, we give the following corollary:
	
	\begin{corollary}
		
		\label{zdgcontent}
		
		Let $S$ be a commutative semiring with a nonzero identity and $M$ a content $S$-semimmodule. If the content function from $M$ to the set of finitely generated ideals of $S$ is onto and the graph $\Gamma_{(\Id(S), c)}(M)$ is non-empty, then it is a connected graph with a diameter at most 3.
		
		\begin{proof}
			
			Let $x,y\in M$ be vertices of the graph $\Gamma_{(\Id(S), c)}(M)$. Since $M$ is a content $S$-semimodule, then $c(x)$ and $c(y)$ are both finitely generated ideals of the semiring $S$ \cite[Proposition 23]{Nasehpour2016}. Clearly, $c(x)c(y)$ is also finitely generated. By assumption, the content function $c$ from $M$ to the set of finitely generated ideals of $S$ is onto. So, there is a $z\in M$ such that $c(z) = c(x)c(y)$. By using Theorem \ref{diamsemigroup}, the proof is complete.\end{proof}
	\end{corollary}
	
	Let us recall that a commutative semigroup $(S,\cdot)$ is called positive ordered if $S$ is a semigroup with the zero 0 and there is a partial order $\leq$ on $S$ such that the following conditions are satisfied:
	
	\begin{enumerate}
		\item The partial order $\leq$ is compatible with the multiplication of the semigroup, i.e. $x \leq y$ implies $xz \leq yz$ for all $x,y,z \in S$,
		\item The partial order is positive, i.e. $0<x$ and $0<y$ imply that $0 < xy$ for all $x,y \in S$.
	\end{enumerate} 
	
	\begin{theorem}
		
		\label{diamsemigroup2}
		
		Let $X$ be a non-empty set, $S$ a positive ordered commutative semigroup with zero, and $f$ a function from $X$ to $S$ with this property that for all $w,z \in X$, if $f(w)f(z)\neq 0$, then there exists a $v \in X$ such that $f(w)f(z) \leq f(v)$, $f(v) \leq f(w)$, and $f(v) \leq f(z)$. Then the graph $\Gamma_{(S,f)}(X)$ is connected with $\diam(\Gamma_{(S,f)}(X)) \leq 3$.

		\begin{proof}
			Let $x,y$ be two distinct vertices of $\Gamma_{(S,f)}(X)$. Therefore, there exists $z,w \in X$ such that $f(z)\neq 0$, $f(w)\neq 0$ and $f(x)f(z)=0$ and $f(y)f(w)=0$. Note that $f(x)\neq 0$ and $f(y)\neq 0$. Now we show that $d(x,y)\leq 3$.
			
			The argument for the case $f(x)f(y)=0$ and the case $f(x)f(y)\neq 0$ while $f(w)f(z)=0$ is the same as the argument in the proof of Theorem \ref{diamsemigroup} and therefore, $d(x,y)\leq 3$.
			
			Now imagine $f(x)f(y)\neq 0$ and $f(z)f(w)\neq 0$. Since by assumption, there exists a $v \in X$ such that $f(z)f(w) \leq f(v)$, $f(v) \leq f(z)$, and $f(v) \leq f(w)$, we have $f(x)f(v)=0$ and $f(v)f(y)=0$ and therefore, $d(x,y)\leq 2$ and the proof is complete.
		\end{proof}

	\end{theorem}
	
	Let us recall that if $M$ is an $S$-semimodule and $N$ is an $S$-subsemimodule of $M$, $[N: M]$ is defined to be the set of all elements $s$ of the semiring $S$ such that $sM\subseteq N$. The proof of the following proposition is straightforward, but we bring it here only for the sake of reference.  
	
	\begin{proposition}
		
		\label{quotientideal}
		Let $S$ be a commutative semiring with a nonzero identity and $M$ an $S$-semimodule. Then the following statements hold:

		\begin{enumerate}
			\item If $N$ is an $S$-subsemimodule of $M$, then $[N:M]$ is an ideal of $S$,
			\item If $P$ and $Q$ are $S$-subsemimodules of the $S$-semimodule $M$, then $$[P:M]\cdot [Q:M] \subseteq [P \cap Q : M],$$
			
			\item If  $P$ and $Q$ are $S$-subsemimodules of the $S$-semimodule $M$ and $P \subseteq Q$, then $$[P:M] \subseteq [Q:M].$$
		\end{enumerate}
		
	\end{proposition}
	
	\begin{corollary}
		
		\label{diamsemimodule}
		Let $S$ be a commutative semiring with a nonzero identity and $M$ be a unital $S$-semimodule. Assume that $q$ is a function from $\Sub(M)$ to $\Id(R)$ with $q(N) = [N:M]$. If the graph $\Gamma_{(\Id(R),q)}(\Sub(M))$ is non-empty, then it is a connected graph whose dimeter is at most 3.
		
		\begin{proof} Use Theorem \ref{diamsemigroup2} and Proposition  \ref{quotientideal}.
		\end{proof}
	\end{corollary}
	
	Let us recall that if $S$ is an idempotent commutative semigroup, then $S$ can be ordered by the following order: $x \leq y$ if $xy = x$ for all $x,y \in S$. Additionally, if $(S,\cdot, 0,1)$ is a monoid with the absorbing element 0, then $S$ is called a bounded semilattice \cite{ChajdaHalasKuhr2007}.

	\begin{proposition} 
		
		Let $(S,\cdot, 0,1)$ be a bounded semilattice and $d$ be the largest element of the poset $S-\{0,1\}$ such that $d^2 = 0$. If $f$ is a function from a set $X$ to $S$ such that the graph $\Gamma_{(S,f)}(X)$ has at least two vertices, then $\diam(\Gamma_{(S,f)}(X))=1$.
		
		\begin{proof}
			Let $x,y$ be vertices of the graph $\Gamma_{(S,f)}(X)$. It is clear that $f(x)$ and $f(y)$ are both nonzero and there are two elements $w,z\in X$ such that $f(x)f(w) = 0$ and $f(y)f(z) =0$. Clearly, these two imply that $f(x) \neq 1$ and $f(y) \neq 1$. Therefore, $f(x) \leq d$ and $f(y) \leq d$, because $d$ is the largest element the poset $S-\{0,1\}$. On the other hand, $f(x)f(y) \leq d^2 =0$. Hence, $\{x,y\}$ is an edge of the graph $\Gamma_{(S,f)}(X)$ and the proof is complete.	
		\end{proof}
		
	\end{proposition}
	
	\begin{corollary}
		Let $S$ be a commutative semiring with an identity and $M$ be a unital $S$-semimodule. Also, let $q$ be the function from $\Sub(M)$ to $Id(R)$ with $q(N) = [N:M]$. If $\mathfrak{m}$ is the only maximal ideal of the semiring $S$ such that $\mathfrak{m} ^ 2 =0$ and the graph $\Gamma_{(Id(S),q)}(\Sub(M))$ has at least two vertices, then its diameter is 1.
	\end{corollary} 
	
	\section{Cycles and Cores of Zero-Divisor Graphs and Their Generalizations}\label{sec:core}
	
	Now we proceed to discuss the cycles of the graph $\Gamma_{(S,f)}(X)$. Let $\Gamma$ be a graph. We denote the set of all vertices adjacent to the vertex $a$ of the graph  $\Gamma$ by $N(a)$. In particular, if $X$ is a non-empty set, $S$ a commutative semigroup with zero, and $f$ a function from $X$ to $S$, then $N(a)$ is the set of all vertices $x\in X-\{a\}$ in the graph $\Gamma_{(S,f)}(X)$ such that $f(x) \neq 0$ and $f(a)f(x)=0$.
	
	\begin{lemma}
		\label{cycle4}
		If $a-x-b$ is a path in a graph $\Gamma$, then either $N(a) \cap N(b) = \{x\}$ or $a-x-b$ is contained in a cycle of the length of at most 4.
		
		\begin{proof}
			Let $a-x-b$ be a path in the graph $\Gamma$. It is obvious that $\{x\} \subseteq N(a) \cap N(b)$. If $N(a) \cap N(b) \neq \{x\}$,
			then there exists a vertex $c$ such that $c\notin \{x,a,b\}$ and $c$ is adjacent to the both vertices $a$ and $b$. So, $a-x-b-c-a$ is a path in $\Gamma$. Hence, $a-x-b$ is contained in a cycle of the length $\leq$ 4.
		\end{proof}
		
	\end{lemma}
	
	\begin{theorem}
		
		\label{uniontrirect}
		
		Let $X$ be a non-empty set, $S$ a commutative semigroup with zero, and $f$ a function from $X$ to $S$. Also, let the graph $\Gamma_{(S,f)}(X)$ have at least three vertices such that for all $a,b,x \in X$ if
		$a-x-b$ is a path in $\Gamma_{(S,f)}(X)$ then $N(a) \cap N(b) \neq \{x\}$. If $\Gamma_{(S,f)}(X)$ is a connected graph with $\diam(\Gamma_{(S,f)}(X)) \leq 3$, then any edge in $\Gamma_{(S,f)}(X)$ is contained in a cycle of the length at most 4 and therefore, $\Gamma_{(S,f)}(X)$ is a union of triangles and rectangles.
				
		\begin{proof}
			Let $a-x$ be an edge in $\Gamma_{(S,f)}(X)$. Since by assumption $\Gamma_{(S,f)}(X)$ is connected with $\diam(\Gamma_{(S,f)}(X)) \leq 3$ and possesses at least three vertices, there exists a vertex $b$ such that either $a-x-b$ or $x-a-b$ is a path in $\Gamma_{(S,f)}(X)$ and in any case, by Lemma \ref{cycle4}, $a-x$ is contained in a cycle of the length of at most 4 and, therefore, is an edge of either a triangle or a rectangle.
		\end{proof}
		
	\end{theorem}
	
	Let us recall that the core of a graph $\Gamma$ is the largest subgraph of $\Gamma$ in which every edge is the edge of a cycle in $\Gamma$ \cite{DeMeyerDeMeyer2005}.
	
	\begin{theorem}
		
		\label{core}
		Let $X$ be a non-empty set, $S$ a commutative semigroup with zero, and $f$ a function from $X$ to $S$. Also, let the graph $\Gamma_{(S,f)}(X)$ have at least three vertices and the function $f$ have this property that for all $x,y\in X$ if $f(x)f(y)\neq 0$ then there exists a $z\in X$ such that $f(z)=f(x)f(y)$. If $\Gamma_{(S,f)}(X)$ contains a cycle, then the core $K$ of $\Gamma_{(S,f)}(X)$ is a union of triangles and rectangles.
		
		\begin{proof}
			Let $a_1 \in K$ and suppose that $a_1$ is a part of neither a triangle nor a rectangle in $\Gamma_{(S,f)}(X)$. So, $a_1$ is a part of a cycle \[C \colon a_1 - a_2 - a_3 - a_4 - \dots - a_n - a_1,\] where $n \geq 5$. Without loss of generality, we can suppose that this is the shortest cycle containing $a_1$ and it follows that $\{a_2,a_4\}$ is not an edge of the graph $\Gamma_{(S,f)}(X)$ and by the definition of the graph $\Gamma_{(S,f)}(X)$, $f(a_2)\cdot f(a_4) \neq 0 $. So, by assumption, there exist a $z \in X$ such that $f(z)=f(a_2)\cdot f(a_4)$. Note that $f(a_1)\cdot f(a_2)=f(a_2)\cdot f(a_3)=0$, so $f(a_1)\cdot f(z)=f(z)\cdot f(a_3)=0$. Therefore, $a_1 - z - a_3$ is a path in $\Gamma_{(S,f)}(X)$. Since $C$ is the shortest cycle of the graph $\Gamma_{(S,f)}(X)$ containing $a_1$, $z=a_2$ and we have $f(a_2)=f(a_2)\cdot f(a_4)$. Now consider $0=f(a_2)\cdot ((fa_4)\cdot f(a_5)) = ((f(a_2)\cdot f(a_4))\cdot f(a_5) = f(a_2)\cdot f(a_5) \neq 0 $, a contradiction. This completes the proof.
		\end{proof}
		
	\end{theorem}

\begin{remark}
	Note that Theorem \ref{core} is related to Theorem 1.5 in \cite{DeMeyerMcKenzieSchneider2002}.
\end{remark}
	
	\subsection*{Acknowledgments}
	
	The author is supported by the Department of Engineering Science at the Golpayegan University of Technology and his special thanks go to the Department for providing all necessary facilities available to him for successfully conducting this research.
	
	\bibliographystyle{plain}

\begin{thebibliography}{15.}
		
		\bibitem{AkbariTavallaeeGhezelahmad2012} Akbari, S., Tavallaee, H.A., Khalashi Ghezelahmad, S.: {\em Intersection graph of submodules of a module},  J. Algebra. Appl,  Vol. 11, No. 1 (2012) 1250019 (8 pages).
		
		\bibitem{AndersonAxtellStickles2011} Anderson, D.F., Axtell, M.C., Stickles, J.A. Jr.: {\em Zero-divisor graphs in commutative rings}, In: Fontana, M., Kabbaj, S.-E., Olberding, B., Swanson, I., eds. Commutative Algebra, Noetherian and Non-Noetherian Perspectives, Springer-Verlag, New York, 2011, 23--45.
		
		\bibitem{AndersonBadawi2017} Anderson D.F., Badawi A.: {\em The zero-divisor graph of a commutative semigroup: a survey}, In: Droste M., Fuchs L., Goldsmith B., Str\''{u}ngmann L. (eds) Groups, Modules, and Model Theory - Surveys and Recent Developments, Springer, Cham, 2017, 23--39.
		
		\bibitem{AndersonLivingston1999} Anderson, D.F., Livingston, P.S.: {\em The zero-divisor graph of a commutative ring}, J. Algebra \textbf{217} (2) (1999), 434--447.
		
		
		\bibitem{AndersonMulay2007} Anderson, D.F., Mulay, S.B.: {\em On the diameter and girth of a zero-divisor graph}, J. Pure Appl. Algebra \textbf{210} (2) (2007), 543--550.
		
		\bibitem{Aubert1953} Aubert, K.E.: {\em On the ideal theory of commutative semi-groups}, Math. Scad., \textbf{1} (1953), 39--54.
		
		\bibitem{AykacAkgunesCevik2019} Ayka\c{c}, S., Akg\"{u}ne\c{s}, N., \c{C}evik, A.S.: {\it Analysis of Zagreb indices over zero-divisor graphs of commutative rings}, Asian-Eur. J. Math., {\bf 12}(06) (2019), Article ID: 2040003.
		
		\bibitem{Beck1988} Beck, I.: {\em Coloring of commutative rings}, J. Algebra \textbf{116}(1) (1988), 208--226.
		
		\bibitem{Bosak1964}  Bos\'{a}k, J.: {\it The graphs of semigroups}, in Theory of Graphs and Its Applications, Proc. Sympos. Smolenice (June 1963), Academic Press, New York, 1965, pp. 119--125.
		
		\bibitem{CannonNeuerburgRedmond2005} Cannon G.A., Neuerburg K.M., Redmond S.P.: {\em  Zero-divisor graphs of nearrings and semigroups}, In: Kiechle H., Kreuzer A., Thomsen M.J. (eds) Nearrings and Nearfields, Springer, Dordrecht, 2005.
		
		\bibitem{ChajdaHalasKuhr2007} Chajda, I., Hala\v{s}, R., K\"{u}hr, J.: {\em Semilattice Structures}, Research and Exposition in Mathematics, 30. Heldermann Verlag, Lemgo, 2007.
		
		\bibitem{ChakrabartyGhoshMukherjeeSen2009} Chakrabarty, I., Ghosh, S., Mukherjee T.K., Sen, M.K.: {\it Intersection graphs of ideals of rings}, Discrete Math 309 (2009) 538--5392.
		
		\bibitem{CoykendallWagstaffSheppardsonSpiroff2012} Coykendall, J., Sather-Wagstaff, S.; Sheppardson, L.; Spiroff, S.: {\em On zero divisor graphs}, Francisco, Christopher (ed.) et al., Progress in commutative algebra 2. Closures, finiteness and factorization. Berlin: Walter de Gruyter, De Gruyter Proceedings in Mathematics, 241--299 (2012). 
		
		\bibitem{CozzensMoazzamiStueckle1995} Cozzens, M., Moazzami, D., Stueckle, S.: The Tenacity of a Graph. In: Proc. Seventh International Conference on the Theory and Applications of Graphs, pp. 1111--1122. Wiley, New York (1995).
		
		\bibitem{CsakanyPollak1969} Cs$\acute{a}$k$\acute{a}$ny, B., Poll$\acute{a}$k, G.: {\it The graph of subgroups of a finite group}, Czechoslovak Math. J. 19 (1969) 241--247.
		
		\bibitem{DeMeyerDeMeyer2005} DeMeyer, F.R., DeMeyer, L.: {\em Zero-divisor graphs of semigroups}, Journal of Algebra, Vol. {\bf 283} (2005), 190--198.
		
		\bibitem{DeMeyerMcKenzieSchneider2002} DeMeyer, F.R., McKenzie, T., Schneider, K.: {\em The Zero-divisor graph of a commutative semigroup}, Semigroup Forum, Vol. {\bf 65} (2002) 206--214.
		
		\bibitem{Diestel2017} Diestel, R.: {\it Graph Theory}, 5th end., Springer-Verlag, Berlin, 2017.
		
		\bibitem{ErdosGoodmanPosa1966} Erd\"{o}s, P., Goodman, A.W., and P\'{o}sa, L.: {\em The representation of a graph by set intersections}, Canad. J. Math. \textbf{18} (1966), 106--112.
		
		\bibitem{EpsteinNasehpour2013} Epstein, N., Nasehpour, P.: {\em Zero-divisor graphs of nilpotent-free semigroups}, J. Algebr Comb, {\bf 37}(3) (2013), 523--543.
		
		\bibitem{Golan1999} Golan, J.S.: {\em Semirings and Their Applications}, Kluwer Academic Publishers, Dordrecht, 1999.
		
		\bibitem{HuckabaKeller1979} Huckaba, J.A., Keller, J.M.: {\it Annihilation of ideals in commutative rings}, Pac. J. Math., Vol. {\bf 83}(2) (1979), 375--379.
		
		\bibitem{ImrichKlavzar2000} Imrich, W., Klav\u{z}ar, S.: {\it Product Graphs: Structure and Recognition}, Wiley, 2000.
		
		\bibitem{Lovasz1971} Lov\'{a}sz, L.: {\em On the cancellation law among finite relational structures}, Period Math Hung {\bf 1} (1971), 145--156.
		
		\bibitem{Lucas2006} Lucas, T.G.: {\em The diameter of a zero-divisor graph}, J. Algebra {\bf 301} (2006), 174--193.
		
		\bibitem{MalakootiRadNasehpour2017} Malakooti Rad, P., Nasehpour, P.: {\em On graphs of bounded semilattices}, preprint, arXiv:1711.01308, 2017, to appear in Math. Notes. 
		
		\bibitem{Nasehpour2021} Nasehpour, P.: {\em Eversible and reversible semigroups and semirings}, Asian-Eur. J. Math., {\bf 14}(1) (2021), Article ID 2150006 (17 pages).
		
		\bibitem{Nasehpour2016} Nasehpour, P.: {\em On the content of polynomials over semirings and its applications}, J. Algebra Appl., {\bf 15}(5) (2016), Article ID 1650088 (32 pages).
		
		\bibitem{Nasehpour2010} Nasehpour, P.: {\em Zero-divisors of content algebras}, Arch. Math. (Brno), {\bf 46} (4) (2010), 237--249.
		
		\bibitem{Osba2016} Osba, E. A.: {\em The intersection graph of finite commutative principal ideal rings}, Acta Math. Acad. Paedagog. Nyh\'{a}zi. (N.S.), {\bf 32} (1) (2016), 15--22.
		
		\bibitem{Redmond2002} Redmond, S.P.: {\em The zero-divisor graph of a non-commutative ring}, Internat. J. Commutative Rings, {\bf 1}(4), (2002) 203--211.
		
		\bibitem{Shen2010} Shen, E. R.: {\it Intersection graphs of subgroups of finite groups}, Czech. Math. J., 60 (135) (2010) 945--950.
		
		\bibitem{Talwar1995} Talwar, S.: {\em Morita equivalence for semigroups}, J. Austral. Math. Soc. Ser. A, Vol. {\bf 59}(1) (1995), 81--111.
		
		\bibitem{Zelinka1975}  Zelinka, B.: {\it Intersection graphs of finite abelian groups}, Czech. Math. J.  25(2) (1975) 171--174.
		
		\bibitem{Zelinka1973}  Zelinka, B.: {\it Intersection graphs of lattices}, Matematick\'{y} \v{c}asopis, {\bf 23}(3) (1973), 216--222.
		
	\end{thebibliography}

\end{document}